\documentclass[leqno]{article}
\usepackage{amssymb,amsmath,amsfonts,amsthm}

 \swapnumbers

\theoremstyle{plain}
  \newtheorem{theo}{Theorem}[section]
  \newtheorem{prop}[theo]{Proposition}
  \newtheorem{lem}[theo]{Lemma}
\theoremstyle{definition}
  \newtheorem{define}[theo]{Definition}
 \newtheorem{remark}[theo]{Remark} 

\newcommand{\into}{\hookrightarrow}
\newcommand{\onto}{\to\mskip-14mu\to} 
 
\renewcommand{\hat}{\widehat}

 \DeclareMathOperator{\alg}{alg}
 \DeclareMathOperator{\et}{et}
 \DeclareMathOperator{\Aut}{Aut}
 
 \DeclareMathOperator{\Ext}{Ext}
 \DeclareMathOperator{\Hom}{Hom}
 
 \DeclareMathOperator{\red}{red}
 
 \DeclareMathOperator{\im}{im}
 \def\del{\delta}
  \def\om{\omega}          \def\Om{\Omega}
     
  \def\gam{\gamma}         
  \def\sig{\sigma}
  \def\g{\mathfrak{g}}    
  \def\a{\mathfrak{a}}
  \def\z{\mathfrak{z}} 
  \def\bc{\mathbb{C}}
 \def\fg{\mathfrak{g}}    \def\fa{\mathfrak{a}}
  \def\simto{\overset{\sim}{\longrightarrow}}
\def\prep{\,^{^{\backslash}}\negmedspace}  %pre-prime   
  \def\prepp{\,^{^{\backslash\!\backslash}}\negmedspace}%double pre-prime
 
\def\beqnn{\begin{equation*}}
  \def\eeqnn{\end{equation*}}
  \def\beqn{\begin{equation}}
  \def\eeqn{\end{equation}}

\begin{document}

  \title{Extensions of Algebraic Groups}

  \author{Shrawan Kumar and Karl-Hermann Neeb}
  \date{}
  \maketitle

\section*{Introduction}

 Let $G$ be a connected complex algebraic group and $A$ an abelian
connected algebraic group, together with an algebraic action of $G$ on $A$ via group
automorphisms. The aim of this note is to study the set of isomorphism classes $ \Ext_{\alg}(G,A)$ of extensions of $G$ by $A$ in the algebraic group category. The following is our main result (cf. Theorem 1.8).

 \begin{theo} For  $G$ and $A$ as above, there exists an exact
sequence of abelian groups:
   \[ 0 \to \Hom (\pi_1([G,G]), A) \to \Ext_{\alg}(G,A)  
\overset{\pi}{\longrightarrow} H^2(\g ,\g_{\red}, \a_u) \to 0\,,
     \] where 
$A_u$ is the unipotent radical of $A$,
$G_{\red}$ is a Levi  subgroup of $G$,   $
\g_{\red},   \g, \a_u$ are the Lie algebras of $G_{\red}, G, A_u$ respectively, and  $H^*(\fg ,\fg_{\red} ,\fa_u )$ is the Lie algebra 
cohomology of the pair $(\fg ,\fg_{\red} )$ with coefficients in the 
$\fg$-module $\fa_u$. 
  \end{theo}

Our next main result is the following analogue of the 
 Van-Est Theorem for the 
algebraic group cohomology (cf. Theorem 2.2).

  \begin{theo}  Let $G$ be a connected  algebraic group and let $\fa$ be a 
finite-dimensional algebraic $G$-module.  Then, for any $p\geq 
0$,
  \[
H^p_{\alg} (G,\fa ) \simeq H^p (\fg ,\fg_{\red} ,\fa ).
  \]
  \end{theo}

 This work was done while the authors were visiting the Fields Institute, Toronto (Canada) in July, 2003, hospitality of which is gratefully acknowledged. 
The first author was partially supported from NSF.

\vskip7ex

By an algebraic group $G$ we mean an affine algebraic group over the field
of complex numbers $\bc$ and the varieties are considered over
$\bc$. The Lie algebra of $G$ is denoted by $L(G)$.

\maketitle

\section{Extensions of Algebraic Groups}
  \begin{define} Let $G$ be an algebraic group and $A$ an abelian
algebraic group, together with an algebraic action of $G$ on $A$ via group
automorphisms, i.e., a morphism  of varieties  $\rho : G\times A\to
A$ such that the induced map $G\to \Aut A$ is a group homomorphism.
Such an $A$ is called an {\it algebraic group with $G$-action}.

By $\Ext_{\alg} (G,A)$ we mean the set of isomorphism classes of
extensions of $G$ by $A$ in the algebraic group category, i.e., 
quotient morphisms $q : \hat G \to G$ with kernel isomorphic to 
$A$ as an algebraic group with $G$-action. 
We obtain on $\Ext_{\alg}(G,A)$ the structure of an abelian group 
by assigning to two extensions $q_i : \hat G_i \to G$ of $G$ by $A$ 
the fiber product extension $\hat G_1 \times_G \hat G_2$ 
of $G$ by $A \times A$ and then applying the group morphism $m_A : A \times A \to A$ 
fiberwise 
to obtain an $A$-extension of~$G$ (this is the Baer sum of two extensions). 
Then $\Ext_{\rm alg}$ assigns to a pair of an algebraic group $G$ 
and an abelian algebraic group $A$ with $G$-action, an abelian group, and this assignment is 
contravariant in $G$ (via pulling back the action of $G$ and the extension) 
and if $G$ is fixed, $\Ext_{\rm alg}(G,\cdot)$ is  a covariant functor from 
the category of abelian algebraic groups with $G$-actions to the category of abelian groups. 
Here we assign to a  $G$-equivariant  morphism $\gamma : A_1 \to A_2$ 
of abelian algebraic groups and an extension 
$q : \hat G \to G$ of $G$ by $A_1$ the extension 
$$ \gamma_*\hat G := (A_2 \rtimes \hat G)/\Gamma(\gamma) \to G, \quad 
[(a,g)] \mapsto q(g), $$
where $\Gamma(\gamma)$ is the graph of $\gamma$ in $A_2 \times A_1$
and the semidirect product refers to the action of $\hat G$ on $A_2$ obtained 
by pulling back the action of $G$ on $A_2$ to $\hat G$. In 
view of the equivariance of $\gamma$, its graph is a normal algebraic subgroup of 
$A_2 \rtimes \hat G$, so that we can form the quotient $\gamma_*\hat G$. 

We define a map
  \[ D: \Ext_{\alg}(G,A) \to \Ext(L(G), L(A)) 
   \] by assigning to an extension
  \[ 1 \to A  \overset{i}{\longrightarrow}\hat{G} \overset{q}{\longrightarrow} G \to 1
  \] of algebraic groups the corresponding extension 
  \[ 0 \to L(A)  \overset{di}{\longrightarrow} L(\hat G) 
 \overset{dq}{\longrightarrow} L(G) \to 0
  \] of Lie algebras. Since $i$ is injective, $di$ is injective. Similarly, $dq$ is surjective. Moreover,  
$\dim G = \dim L(G)$  and hence the above
sequence of Lie algebras is indeed exact. 
\end{define}

It is clear from the definition of $D$ that it is a homomorphism of abelian 
groups. If $\g$ is the Lie algebra of $G$ and $\a$ the Lie algebra of $A$, then 
the group $\Ext(\g,\a)$ is isomorphic to the second Lie algebra cohomology space 
$H^2(\g,\a)$ of $\g$ with coefficients in the $\g$-module $\a$ (with respect to the derived 
action) ([CE]). 
Therefore the description of the group $\Ext_{\rm alg}(G,A)$ depends on a 
good description of kernel and cokernel of $D$ which will be obtained below in 
terms of an exact sequence involving $D$. 

In the following $G$ is always assumed to be connected. 
The following lemma reduces the extension theory for connected algebraic groups $A$ with $G$-actions
to the two cases of a torus $A_s$ and the case of a unipotent group $A_u$. 

\begin{lem} Let $G$ be connected and $A$ be a connected algebraic group with $G$-action. 
Further, let $A = A_u A_s$ denote the decomposition of $A$ into its
unipotent and  reductive factors. Then $A \cong A_u \times A_s$ as a $G$-module, where 
$G$ acts trivially on $A_s$ and $G$ acts on $A_u$ as a $G$-stable subgroup
of $A$. Thus,
we have 
 \begin{equation} \Ext_{\rm alg}(G,A) \cong   \Ext_{\rm alg}(G,A_u) \oplus \Ext_{\rm alg}(G,A_s). 
\end{equation} 
\end{lem} 

\begin{proof} Decompose
  \begin{equation} A = A_uA_s ,
  \end{equation} where $A_s$ is the set of semisimple elements of $A$ and
$A_u$ is the set of unipotent elements of $A$.  Then $A_s$ and $A_u$ are
closed subgroups of $A$ and (2) is a direct product decomposition  (see
[H, Theorem 15.5]).  The action of $G$ on
$A$ clearly keeps $A_s$ and $A_u$ stable separately.  Also, $G$ acts
trivially on $A_s$ since $\Aut (A_s)$ is discrete and $G$ is connected (by assumption).  Thus the action of $G$
on $A$ decomposes as the product of actions on $A_s$ and $A_u$ with the trivial action on 
$A_s$.  Hence the isomorphism (1)
 follows from the functoriality of $\Ext_{\rm alg}(G,\cdot)$. 
\end{proof} 

If $G = G_u \rtimes G_{\rm red}$ is a Levi decomposition of $G$, then
$G_u$ being simply-connected, 
$$\pi_1(G) \cong \pi_1(G_{\rm red}), $$
where $G_u$ is the unipotent radical of $G$, $G_{\rm red}$ is a Levi subgroup of $G$ and $\pi_1$ denotes the fundamental 
group.
The connected reductive group $G_{\rm red}$ is a product of its connected center 
$Z := Z(G_{\rm red})_0$ and its commutator group 
$G_{\rm red}' := [G_{\rm red}, G_{\rm red}]$ which is a connected semisimple group. 
Thus,  $G_{\rm red}'$ has an algebraic universal covering group 
$\tilde G_{\rm red}'$, with
the finite abelian group $\pi_1(G_{\rm red}')$ as its fiber. 
We write $\tilde G_{\rm red} := Z \times \tilde G_{\rm red}'$ which is  
 an algebraic covering group of $G_{\rm red}$; denote its kernel 
by $\Pi_G$ and observe that 
$$ \tilde G := G_u \rtimes \tilde G_{\rm red} $$
is a covering of $G$ with $\Pi_G$ as its fiber. We write $q_G  :  \tilde G \to G$ for the corresponding 
covering map. 

\begin{lem} If $G$ and $A$ are tori, then $\Ext_{\rm alg} (G,A) = 0$. 
\end{lem} 

\begin{proof} Let $q : \hat G \to G$ be an extension of the torus 
$G$ by $A$. Then, as is well known, $\hat G$ is again a torus (cf. [B, $\S$11.5]).  Since 
any character of a subtorus of a torus extends to a character of the whole 
groups ([B, \S 8.2]),  
the identity $I_A : A \to A$ extends to a morphism 
$f : \hat G \to A$.  Now  $\ker f$ yields a splitting of the above extension.   
\end{proof}

The following proposition deals with the case $A = A_s$. 

\begin{prop} If $A = A_s$, then $D = 0$ and we obtain an exact sequence 
$$ \Hom(\tilde G, A_s) \overset{res}{\longrightarrow}\Hom(\Pi_G, A_s) 
\overset{\Phi}{\longrightarrow} \Ext_{\rm alg}(G, A_s), $$
where $\Phi$ assigns to any $\gamma \in \Hom(\Pi_G, A_s)$ the extension $\gamma_* \tilde G$. 
The kernel of $\Phi$ consists of  those homomorphisms vanishing on 
the fundamental group  $\pi_1(G_{\rm red}')$ of $G_{\rm red}'$ and $\Phi$ factors through an isomorphism 
$$ \Phi'  :  \Hom(\pi_1(G_{\rm red}'), A_s) \simeq \Ext_{\rm alg} (G, A_s). $$
\end{prop} 

\begin{proof} Consider an extension  
\[1\to A_s\to \hat{G} \to
G\to 1.\]
Since $A_s$ is a central torus in $\hat G$, the unipotent radical $\hat{G}_u$ of
$\hat{G}$ maps isomorphically on $G_u$.  Also 
\[1\to A_s \to \hat{G}_{\red}
\to G_{\red} \to 1\]
is an extension whose restriction to $Z$ splits by the preceding lemma. 
On the other hand the commutator group of $\hat G_{\rm red}$ has the same 
Lie algebra as $G_{\rm red}'$, hence is a quotient of $\tilde G_{\rm red}'$. 
Thus $\hat G_{\rm red}$ is a quotient of 
$A_s \times Z \times \tilde G_{\rm red}'$, which implies 
that $\hat G$ is a quotient of $A_s \times \tilde G$. Hence $\hat G$ is obtained from  $A_s \times \tilde G$ via taking its quotient by 
the graph of a homomorphism $\Pi_G \to A_s$. Conversely, any such extension $\hat G$ of $G$ is obtained this way. This proves that $\Phi$ is surjective.
In particular,  
 the pullback $q_G^*\hat G$ of $\hat G$ 
to $\tilde G$ always splits. 

We next show  that $\ker\Phi$ coincides with the image of the restriction map from 
$\Hom(\tilde G, A_s)$ 
to $\Hom(\Pi_G, A_s)$.  Assume that the extension $\hat G_\gamma = \gamma_*\tilde G$ 
defined by 
$\gamma \in \Hom(\Pi_G, A_s)$ splits. 
Let $\sigma : G \to \hat G_\gamma$ be a splitting morphism. 
Pulling $\sigma$ back via $q_G$, we obtain a splitting morphism 
$$\tilde \sigma  :  \tilde G \to q_G^*\hat G_\gamma \cong A_s \times \tilde G.  $$
 Thus, there exists a morphism $\delta  :  \tilde G \to A_s$ 
of algebraic groups such that  $\sigma$ satisfies 
$\sigma(q_G(g)) = 
\beta(\delta(g),g)$  for all $g\in \tilde G$, where $\beta:
 A_s \times \tilde G \to  \hat G_\gamma =  (A_s \times \tilde G)/\Gamma (\gamma)$ is the standard quotient map. For $g \in \Pi_G = \ker q_G$ we have 
$\beta(\delta(g), g) = 1$, and therefore $\delta(g) = \gamma(g)$ 
for all $g\in \Pi_G$.
 This shows that $\delta$ is an extension of $\gamma$ to $\tilde G$. 
Conversely, if $\gamma$ extends to $\tilde G$, $\hat G_\gamma$ is a trivial extension of $G$. 

That $D = 0$ follows from the fact that 
$\hat G$ and $q_G^* \hat G$ have the same Lie algebras, which is a 
split extension of $\g$ by $\a_s$. 

We recall that $\tilde G = G_u \rtimes (Z \times \tilde G_{\rm red}')$. 
If a homomorphism $\gamma  :  \Pi_G \to A_s$ extends to 
$\tilde G$, then it must vanish on the subgroup 
$\pi_1(G_{\rm red}')$ of $\Pi_G$ since, $\tilde G_{\rm red}'$ being a semisimple group, there are no nonconstant homomorphisms from
$\tilde G_{\rm red}' \to A_s$. Conversely, if a homomorphism
$\gamma  :  \Pi_G \to A_s$ vanishes on $\pi_1(G_{\rm red}')$, 
 then $\gamma$ defines a homomorphism 
$$ Z \cap G_{\rm red}' \cong \Pi_G/ \pi_1(G_{\rm red}') \to A_s. $$
But $A_s$ being a torus, 
this extends to a morphism $f  :  Z \to A_s$ ([B, \S 8.2]) which in turn can be pulled 
back via $Z \cong \tilde G/(G_u \rtimes \tilde G_{\rm red}')$ to a morphism 
$\tilde f  :  \tilde G \to A_s$ extending $\gamma$. This proves that 
the image of $\Hom(\tilde G, A_s)$ 
under the restriction map in $\Hom(\Pi_G, A_s)$ is 
the annihilator of $\pi_1(G_{\rm red}')$, so that 
$$ \Phi  :  \Hom(\Pi_G, A_s) \to \Ext_{\rm alg}(G, A_s) $$
factors through an isomorphism 
$$ \Phi'  :  \Hom(\pi_1(G_{\rm red}'), A_s) \simeq \Ext_{\rm alg}(G, A_s). $$
\end{proof} 

\begin{remark} A unipotent group $A_u$ over $\bc$  has no 
non-trivial finite subgroups, 
so that 
$$ \Hom(\pi_1(G_{\rm red}'), A_s)
\cong \Hom(\pi_1(G_{\rm red}'), A). $$
\end{remark}

Now we turn to the study of extensions by unipotent groups. In contrast to the situation 
for tori, we shall see that these extensions are faithfully represented 
by the corresponding Lie algebra extensions. 

\begin{lem} The canonical restriction map
\[ H^2(\g ,\g_{\red}, \a_u) \longrightarrow H^2(\g ,\a_u)
\] is injective.
\end{lem}

\begin{proof} Let $\omega \in Z^2(\g,\a_u)$ be a Lie algebra cocycle 
 representing an element of 
$H^2(\g,\g_{\red},\a_u)$ and suppose that the class $[\omega] \in H^2(\g,\a_u)$ vanishes, 
so that the extension 
$$ \hat\g := \a_u \oplus_\omega \g \to \g, \quad (a,x) \mapsto x $$
with the bracket $[(a,x), (a',x')] = (x.a' - x'.a + \omega(x,x'), [x,x'])$ 
splits. We have to find a $\g_{\red}$-module  map $f  :  \g \to \a_u$ vanishing on $\g_{\red}$ with 
$$\omega(x,x') = (d_\g f)(x,x') := x.f(x') -x'.f(x) - f([x,x']), \quad x,x' \in \g. $$
Since the space $C^1(\g,\a_u)$  of linear maps $\g \to \a_u$ 
is  a semisimple $\g_{\red}$-module ($\a_u$ being a $G$-module, in particular, a $G_{\red}$-module), we have 
$$ C^1(\g,\a_u) = C^1(\g,\a_u)^{\g_{\red}} \oplus \g_{\red}.C^1(\g,\a_u) $$
and similarly for the space $Z^2(\g,\a_u)$ of $2$-cocycles. 
As the Lie algebra differential $d_\g  :  C^1(\g,\a_u) \to Z^2(\g, \a_u)$ 
is a $\g_{\red}$-module map, each $\g_{\red}$-invariant coboundary is the image of a 
$\g_{\red}$-invariant cochain in $C^1(\g,\a_u)$. 
We conclude, in particular, that $\omega = d_\g h$ for some 
$\g_{\red}$-module map  $h : \g \to \a_u$. 
For $x \in \g_{\red}$ and $x' \in \g$ it follows that 
\begin{eqnarray*}
0 &=& \omega(x,x') 
= x.h(x') - x'.h(x) - h([x,x']) \\ 
&=& h([x,x']) - x'.h(x) - h([x,x']) = - x'.h(x), \\ 
\end{eqnarray*}
showing that $h(\g_{\red}) \subseteq \a_u^\g$, which in turn leads to 
$[\g_{\red}, \g_{\red}] \subseteq \ker h$. 
As $\z(\g_{\red}) \cap [\g,\g] = \{0\}$, 
the map $h \vert_{\z(\g_{\red})}$ extends to a linear map 
$f  :  \g \to \a_u^\g$ vanishing on $[\g,\g]$. Moreover, since $f$ vanishes on $[\g,\g] $, $f$ is clearly a $\g$-module map, in particular, a 
$\g_{\red}$-module map. Then $d_\g f = 0$, so that 
$d_\g(h - f) = \omega$, and $h - f$ vanishes on $\g_{\red}$. 
 \end{proof}

\begin{prop} For $A = A_u$ the map $D:\Ext_{\rm alg}(G, A_u) \to 
H^2(\g,  \a_u)$ induces a bijection 
$$ D  :  \Ext_{\rm alg}(G, A_u) \to H^2(\g, \g_{\rm red}, \a_u). $$
\end{prop}

\begin{proof} In view of the preceding lemma, we may identify 
$H^2(\g,\g_{\red},\a_u)$ with a subspace of $H^2(\g,\a_u)$. 
First we claim that $\im(D)$ is contained in this subspace. 
For any extension
  \begin{equation} 1 \to A_u \to \hat{G} \to G \to 1, 
  \end{equation} we choose a Levi subgroup $\hat{G}_{\red}\subset
\hat{G}$ mapping to $G_{\rm red}$ under the above map $\hat{G} \to G$.  Then
  \[ \hat{G}_{\red} \cap A_u = \{1\}.
  \] Moreover, $\hat{G}_{\red} \to G_{\red}$ is surjective and hence an
isomorphism.  This shows that the extension (3) restricted to $G_{\red}$
is trivial and that $\hat\g_u$ contains a $\hat\g_{\red}$-invariant complement to 
$\a_u$. Therefore $\hat\g$ can be described by a cocycle $\omega \in Z^2(\g,\g_{\red},\a_u)$, in particular, $\omega$ 
vanishes on $\g \times \g_{\red}$. This shows that Im $D \subset 
H^2(\g, \g_{\rm red}, \a_u)$.

If the image of the extension (3) under $D$ vanishes, then 
the extension $\a_u \into \hat\g_u \onto \g_u$ splits, which implies that 
the corresponding extension of unipotent groups 
$A_u \into \hat G_u \onto G_u$ splits. Moreover, 
the splitting map can be chosen to be  $G_{\red}$-equivariant, since
$\omega$ is $G_{\red}$-invariant. This means that we have a morphism 
$G_u \rtimes G_{\red} \to \hat G \cong \hat G_u \rtimes G_{\red}$ 
splitting the extension (3). This proves that $D$ is injective. 

To see that $D$ is surjective, let $\omega \in Z^2(\g,\g_{\red},\a_u)$. 
Let $q  : \hat\g := \a_u \oplus_\omega \g \to \g$ denote the corresponding Lie 
algebra extension. Since $\a_u$ is a nilpotent module of $\g_u$, the subalgebra 
$\hat\g_u := \a_u \oplus_\omega \g_u$ of $\hat\g$ is nilpotent, hence corresponds 
to a unipotent algebraic group $\hat G_u$ which is an extension of 
$G_u$ by $A_u$. Further, the $G_{\red}$-invariance of the decomposition 
$\hat\g = \a_u \oplus \g$ implies that $G_{\red}$ acts  algebraically 
on $\hat\g_u$ and hence on $\hat G_u$, so that we can form the semidirect product 
$\hat G := \hat G_u \rtimes G_{\red}$ which is an extension of $G$ by $A_u$ 
mapped by $D$ onto~$\hat\g$. 
\end{proof}

 \begin{theo} For a connected algebraic group $G$ and a connected abelian
algebraic group $A$ with $G$-action, there exists an exact
sequence of abelian groups:
   \[ 0 \to \Hom (\pi_1([G,G]), A) \to \Ext_{\alg}(G,A)  
\overset{\pi}{\longrightarrow} H^2(\g ,\g_{\red}, \a_u) \to 0\,,
     \] where $\a = L(A)$, 
$G_{\red}$ is a Levi  subgroup of $G$,  
$\g_{\red} = L(G_{\red})$,   $\g = L(G)$ and $\a_u = L(A_u)$. 

(Observe that, by the following proof, the fundamental group $\pi_1([G,G])$ is a finite group.)

  \end{theo}
 
   \begin{proof} In view of the Levi decomposition of  the commutator $[G,G] = 
[G,G]_u \rtimes G_{\red}'$, 
we have $\pi_1([G,G]) = \pi_1(G_{\red}')$. 
Now we only have to use Lemma 1.2 to combine the preceding results 
Propositions 1.4 and 1.7 on extensions by 
$A_s$ and $A_u$ to complete the proof. 
\end{proof}

\section{Analogue of Van-Est Theorem for algebraic group 
cohomology}

  \begin{define}  Let $G$ be an algebraic group and $A$ an 
abelian algebraic group with $G$-action.
For any $n\geq 0$, let $C^n_{\alg} (G,A)$ be the abelian group 
consisting of all the variety morphisms $f: G^n \to A$ under the
pointwise addition.  Define the differential
  \begin{align*}
\del : C^n_{\alg} (G,A) &\to C^{n+1}_{\alg} (G,A)\qquad \text{by}\\
(\del f) (g_0,\cdots ,g_n) &= g_0\cdot f(g_1,\cdots , g_n) + 
(-1)^{n+1} f(g_0,\cdots ,g_{n-1})\\
& \qquad + \sum^{n-1}_{i=0} (-1)^{i+1}\, f(g_0,g_1,\cdots 
,g_ig_{i+1},\cdots ,g_n) .
  \end{align*}
Then, as is well known (and easy to see), 
  \beqn
\del^2 = 0 .
  \eeqn

The {\em algebraic group cohomology  $H^*_{\alg} (G,A)$ of $G$ 
with coefficients in $A$} is defined as the cohomology of the 
complex
  \[
0 \to C^0_{\alg} (G,A) \overset{\del}{\longrightarrow} C^1_{\alg} 
(G,A) \overset{\del}{\longrightarrow} \cdots .
  \]
  \end{define}

We have the following analogue of the Van-Est Theorem [V] for the 
algebraic group cohomology.

  \begin{theo}  Let $G$ be a connected algebraic group and let $\fa$ be a 
finite-dimensional algebraic $G$-module.  Then, for any $p\geq 
0$,
  \[
H^p_{\alg} (G,\fa ) \simeq H^p (\fg ,\fg_{\red} ,\fa ) ,
  \]
where $\fg$ is the Lie algebra of $G$ and  $\fg_{\red}$ is the Lie 
algebra of a Levi subgroup $G_{\red}$ of $G$ as in 
Section 1.
  \end{theo}

  \begin{proof}  Consider the homogeneous affine variety $X := 
G/G_{\red}$ and let $\Om^q (X,\fa )$ denote the complex vector 
space of algebraic de Rham forms on $X$ with values in the vector 
space $\fa$.  Since $X$ is a $G$-variety under the left 
multiplication of $G$ and $\fa$ is a $G$-module, $\Om^q$ has a 
natural locally-finite algebraic $G$-module structure.  Define a 
double cochain complex $A = \bigoplus_{p,q\geq 0} A^{p,q}$, where
  \beqnn
A^{p,q} := C^p_{\alg} (G, \Om^q(X,\fa ))
  \eeqnn
and $C^p_{\alg} (G,\Om^q(X,\fa ))$ consists of all the maps $f: G^p 
\to \Om^q(X,\fa )$ such that $\im f\subset M_f$, for some 
finite-dimensional $G$-stable subspace $M_f\subset \Om^q(X,\fa )$ 
and, moreover, the map $f: G^p\to M_f$ is algebraic.  Let $\del : 
A^{p,q} \to A^{p+1,q}$ be the group cohomology differential as in 
Section 2.1 and let $d: A^{p,q} \to A^{p,q+1}$ be induced from 
the standard de Rham differential $\Om^q(X,\fa )\to 
\Om^{q+1}(X,\fa )$, which is a $G$-module map.  It is easy to see 
that $d\del -\del d=0$ and, of course, $d^2=\del^2=0$.  Thus, 
$(A,\del ,d)$ is a double cochain complex.  This gives rise to 
two spectral sequences both converging to the cohomology of the 
associated single complex $(C, \del +d)$ with their $E_1$-terms 
given as follows:
  \begin{align}
 \,^{\backslash}E^{p,q}_1 &= H^q_d(A^{p,*}),\qquad\text{and}\\
 {\prepp}E_1^{p,q} &= H^q_\del (A^{*,p}).
\end{align}
We now determine $\prep E_1$ and ${\prepp}E_1$ more explicitly in 
our case.

Since $X$ is a contractible variety, by the algebraic de Rham theorem 
[GH, Chap. 3, \S 5], the algebraic deRham cohomology
   \beqnn
H^q_{\text{dR}}(X,\fa )  \begin{cases}
\simeq \fa , &\text{if $q=0$}\\
= 0, &\text{otherwise}.  
  \end{cases}    \eeqnn
Thus,
   \beqnn
\prep E^{p,q}_1  \begin{cases}
\simeq C^p_{\alg}(G, \fa ), &\text{if $q=0$}\\
= 0, &\text{otherwise}.
  \end{cases}    \eeqnn
Therefore, 
   \beqn
\prep E_2^{p,q}= H^p_{\del}(H^q_d(A))=\begin{cases}
 H^p_{\alg}(G,\fa ), &\text{if $q=0$}\\
 0, &\text{otherwise}.   \end{cases}    
   \eeqn
In particular, the spectral sequence $\prep E_*$ collapses at $\prep E_2$.  
From this we see that there is a canonical isomorphism
  \beqn
H^p_{\alg}(G,\fa ) \simeq H^p(C, \del +d).
  \eeqn

We next determine $\prepp E_1$ and $\prepp E_2$.  But first we need the 
following two lemmas.

  \begin{lem}
For any $p\geq 0$,
  \[
H^q_{\alg}(G,\Om^p (X,\fa )) =\begin{cases}  \Om^p(X,\fa )^G, &\text{if 
$q=0$}\\ 0, &\text{otherwise},  \end{cases}
  \]
where $\Om^p(X,\fa )^G$ denotes the subspace of $G$-invariants in 
$\Om^p(X,\fa )$.
  \end{lem}

  \begin{proof}
The assertion for $q=0$ follows from the general properties of group 
cohomology. So we need to consider the case $q>0$ now.

Since $L := G_{\red}$ is reductive, any algebraic $L$-module $M$ is 
 completely reducible.
Let
  \[
\pi^M : M \to M^L
  \]
be the unique $L$-module projection onto the space of $L$-module 
invariants $M^L$ of $M$.  Taking $M$ to be the ring of regular functions 
$\bc [L]$ on $L$ under the left regular representation, i.e., under the 
action 
  \[
(k\cdot f)(k') = f(k^{-1}k^1), \quad\, \text{for }f\in\bc [L], k,k'\in L,
  \]
we get the $L$-module projection $\pi =\pi^{\bc [L]}: \bc [L] \to \bc$.
Thus, for any complex vextor space $V$, we get the projection
$\pi\otimes I_V : \bc [L]\otimes V\to V$, which we abbreviate simply by
$\pi$, where $I_V$ is the identity map of $V$.
 We define a `homotopy operator' $H$, for any $q\geq 0$, 
  \[
H: C^{q+1}_{\alg}(G,\Om^p(X,\fa )) \to C^q_{\alg}(G, \Om^p(X,\fa ))
  \]
by
  \[
\Bigl( (Hf) (g_1,\cdots,g_q)\Bigr)_{g_0L} = \pi 
\bigl(\Theta^f_{(g_0,\cdots,g_q)}\bigr),
  \]
for $f\in C^{q+1}_{\alg}(G,\Om^p(X,\fa ))$ and $g_0,\cdots ,g_q\in G$, 
where $\Theta^f_{(g_0,\cdots ,g_q)} : L\to \Om^p(X,\fa )_{g_0L}$ is defined 
by 
  \[
\Theta^f_{(g_0,\cdots ,g_q)}(k) = \Bigl( (g_0k)\cdot f(k^{-1}g_0^{-1}, 
g_1,g_2,\cdots ,g_q)\Bigr)_{g_0L},
  \]
for $k\in L$.  (Here $\Om^p(X,\fa )_{g_0L}$ denotes the fiber at $g_0L$ of 
the vector bundle of $p$-forms in $X$ with values in $\fa$ and, for a form 
$\om$, $\om_{g_0L}$ denotes the value of the form $\om$ at $g_0L$.)
It is easy to see that on $C^{q}_{\alg}(G,\Om^p(X,\fa))$, for any $q\geq 
1$, 
  \beqn
H\del +\del H=I. 
  \eeqn
 To prove this, take  any $f\in C^q_{\alg}(G, \Om^p(X,\fa ))$ and 
 $g_0,\cdots ,g_q\in G$. Then, 
\begin{align}
\Bigl( (H\del f) (g_1,\cdots,g_q)\Bigr)_{g_0L} 
&= \pi 
\bigl(\Theta^{\del f}_{(g_0,\cdots,g_q)}\bigr)\notag\\
&=\bigl(f(g_1,\cdots,g_q)\bigr)_{g_0L}\notag\\
&+(-1)^{q+1}
\pi\Bigl( \bigl((g_0k)\cdot f(k^{-1}g_0^{-1}, 
g_1,\cdots ,g_{q-1})\bigr)_{g_0L}\Bigr)\notag\\
&+\sum_{i=1}^{q-1} \, (-1)^{i+1}
\pi\Bigl( \bigl((g_0k)\cdot f(k^{-1}g_0^{-1}, 
g_1, \cdots, g_ig_{i+1}, \cdots, g_{q})\bigr)_{g_0L}\Bigr)\notag\\
&-\pi\Bigl( \bigl((g_0k)\cdot f(k^{-1}g_0^{-1}g_1, 
g_2, \cdots, g_{q})\bigr)_{g_0L}\Bigr),
  \end{align}
where $\bigl((g_0k)\cdot f(k^{-1}g_0^{-1}, 
g_1,\cdots ,g_{q-1})\bigr)_{g_0L}$ means the function from 
$L$ to $\Om^p(X,\fa )_{g_0L}$ defined as $k\mapsto \bigl((g_0k)\cdot f(k^{-1}g_0^{-1}, 
g_1,\cdots ,g_{q-1})\bigr)_{g_0L}$. Similarly,
 \begin{align}
\bigl( (\del Hf) (g_1,\cdots,g_q)\bigr)_{g_0L} 
&=
\Bigl( g_1\cdot \bigl((Hf) (g_2,\cdots,g_q)\bigr)\Bigr)_{g_0L}\notag\\ 
&+(-1)^q \bigl((Hf) (g_1,\cdots, g_{q-1})\bigr)_{g_0L}\notag\\
&+\sum_{i=1}^{q-1} \, (-1)^{i}
\bigl( (Hf) (g_1, \cdots, g_ig_{i+1}, \cdots, g_{q})\bigr)_{g_0L}\notag\\
&=\Bigl( g_1\cdot \bigl((Hf) (g_2,\cdots,g_q)\bigr)\Bigr)_{g_0L}\notag\\
&+ (-1)^{q}
\pi\Bigl( \bigl((g_0k)\cdot f(k^{-1}g_0^{-1}, 
g_1,\cdots ,g_{q-1})\bigr)_{g_0L}\Bigr)\notag\\
&+\sum_{i=1}^{q-1} \, (-1)^{i}
\pi\Bigl( \bigl((g_0k)\cdot f(k^{-1}g_0^{-1}, 
g_1, \cdots, g_ig_{i+1}, \cdots, g_{q})\bigr)_{g_0L}\Bigr).
\end{align}
From the definition of the $G$-action on  $\Om^p(X,\fa )$, it is easy to see that 
\beqn
\pi\Bigl( \bigl((g_0k)\cdot f(k^{-1}g_0^{-1}g_1, 
g_2, \cdots, g_{q})\bigr)_{g_0L}\Bigr)= \Bigl( g_1\cdot \bigl((Hf) (g_2,\cdots,g_q)\bigr)\Bigr)_{g_0L}.
\eeqn
Combining (10)-(12), we clearly get (9).

From the above identity (9), we see, of course,  that any cocycle 
in $C^q_{\alg}(G, \Om^p(X,\fa ))$ (for any $q \geq 1$) is a coboundary, 
proving the lemma.  
  \end{proof}

  \begin{lem}
The restriction map $\gam : \Om^p(X,\fa )^G\to C^p(\fg ,\fg_{\red},\fa )$ 
 (defined below in the proof) is an isomorphism for all $p\geq 0$, where 
$ C^*(\fg ,\fg_{\red},\fa )$ is the standard cochain complex for the Lie algebra pair $(\fg ,\fg_{\red})$ with coefficient in 
the $\fg$-module~$\fa$.  
Moreover, $\gam$ commutes with differentials.  Thus, $\gam$ induces an 
isomorphism in cohomology
  \[
H^*(\Om (X,\fa )^G) \simto H^*(\fg ,\fg_{\red},\fa).
  \]
  \end{lem}

  \begin{proof}
For any $\om\in\Om^p(X,\fa )^G$, define $\gam (\om )$ as the value of 
$\om$ at $eL$.  Since $G$ acts transitively on $X$, and $\om$ is 
$G$-invariant, $\gam$ is injective.

Since any $\om_o\in C^p(\fg ,\fg_{\red},\fa )$ can be extended (uniquely) 
to a $G$-invariant form on $X$ with values in $\fa$, $\gam$ is surjective.  
Further, from the definition of differentials on the two sides, it is easy 
to see that $\gam$ commutes with differentials.
  \end{proof}
  \medskip

\noindent {\bfseries 2.5 Continuation of the proof of Theorem 2.2.}

We now determine $\prepp E$.  First of all, by (6) of (2.2), 
  \begin{align*}
\prepp E^{p,q}_1 &= H^q_{\del}(A^{*,p})
= H^q_{\alg}(G,\Om^p(X,\fa )).
  \end{align*}
Thus, by Lemma 2.3,
  \[
\prepp E_1^{p,0} = H^0_{\alg}(G, \Om^p(X,\fa )) = \Om^p(X,\fa )^G,
  \]
and
  \[
\prepp E_1^{p,q} = 0, \qquad\text{if } q>0.
  \]

Moreover, under the above equality, the  differential of the spectral 
sequence $d_1 : \prepp E^{p,0}_1 \to \prepp E_1^{p+1,0}$ can be identified 
with the restriction of the deRham differential
  \[
\Om^p (X,\fa )^G \to \Om^{p+1}(X,\fa )^G.
  \]

Thus, by Lemma 2.4,
  \beqn  
\prepp E^{p,q}_2 =\begin{cases}  H^p (\fg ,\fg_{\red} ,\fa ), 
&\text{if }q=0 \\  0, &\text{otherwise}.  \end{cases}
  \eeqn

In particular, the spectral sequence $\prepp E$ as well degenerates at the 
$\prepp E_2$-term.  Moreover, we have a canonical isomorphism
  \beqn 
H^p(\fg ,\fg_{\red}, \fa ) \simeq H^p(C, \del +d).
  \eeqn
Comparing the above isomorphism with the isomorphism (8) of \S 2.2, we get 
a canonical isomorphism:
  \[
H^p_{\alg}(G,\fa ) \simeq H^p(\fg ,\fg_{\red}, \fa ).
  \]
This proves Theorem 2.2.
  \end{proof}

\begin{remark} Even though we took the field $\bc$ as our base field, all the results of this paper hold (by the same proofs)
over any algebraically closed field of char. $0$, if we replace the fundamental group $\pi_1$ by the algebraic fundamental group.  
\end{remark}

\vskip5ex

Addresses:  Shrawan Kumar, Department of Mathematics,
University of North Carolina,
Chapel Hill, NC 27599-3250, USA\\
shrawan$@$email.unc.edu\\

Karl-Hermann Neeb, Fachbereich Mathematik,
Darmstadt University of Technology, 
Schlo$\beta$gartenstr. 7,
D-64289, Darmstadt,
Germany\\
 neeb$@$mathematik.tu-darmstadt.de

\end{document}